\documentclass[12pt,a4paper,reqno]{amsart}
 
\usepackage{amsmath,amssymb,amsthm}

\numberwithin{equation}{section}

\newtheorem{theorem}[equation]{Theorem}
\newtheorem{lemma}[equation]{Lemma}
\newtheorem{corollary}[equation]{Corollary}
\newtheorem{proposition}[equation]{Proposition}

\theoremstyle{definition}
\newtheorem{definition}[equation]{Definition}
\newtheorem{example}[equation]{Example}
\newtheorem{remark}[equation]{Remark}

\newcommand{\N}{\mathbb N}

\newcommand{\R}{\mathbb R}
\newcommand{\Z}{\mathbb Z}
\newcommand{\C}{\mathbb C}
\newcommand{\Q}{\mathbb Q}
\newcommand{\T}{\mathbb T}
\newcommand{\fii}{\varphi}
\newcommand{\bddlin}{\mathcal L}
\newcommand{\tr}{{\rm tr}}
\newcommand{\lin}{{\rm lin}}
\newcommand{\range}{\mathsf{R}}
\newcommand{\domain}{\mathsf{D}}

\newcommand{\no}[1]{\| #1 \|}
\newcommand{\Cinfty}{\mathsf{C}^{\infty}}
\newcommand{\rank}{\operatorname{rank}}
\newcommand{\<}{\langle}
\renewcommand{\>}{\rangle}
\newcommand{\ket}[1]{\left|#1\right\rangle}
\newcommand{\bra}[1]{\left\langle#1\right|}
\newcommand{\kb}[2]{\left|#1\right\rangle\left\langle#2\right|}

\title[Generalized eigenvalue expansions]{Positive sesquilinear form measures and generalized eigenvalue expansions}
\author[Hyt\"onen \and Pellonp\"a\"a \and Ylinen]{}
\date{\today}

\begin{document}

\maketitle

\begin{center}

Tuomas Hyt\"onen\\
Department of Mathematics and Statistics\\ University of Helsinki\\ Gustaf H\"allstr\"omin katu 2b\\ FI-00014 Helsinki, Finland\\ \texttt{tuomas.hytonen@helsinki.fi}\\ \mbox{}\\

Juha-Pekka Pellonp\"a\"a\footnote{Corresponding author. Telephone: +358-2-333 5737, telefax: +358-2-333 5070.}\\

Department of Physics\\ University of Turku\\ FI-20014 Turku, Finland\\  \texttt{juhpello@utu.fi}\\ \mbox{}\\

Kari Ylinen\\
Department of Mathematics\\ University of Turku\\ FI-20014 Turku, Finland\\ \texttt{ylinen@utu.fi}

\end{center}

\newpage

\subsection*{Abstract}
Positive operator measures (with values in the space of bounded
operators on a Hilbert space) and their generalizations, mainly
positive sesquilinear form measures,  are considered with the aim
of providing a framework for their generalized eigenvalue type
expansions. Though there are formal similarities with earlier
approaches to special cases of the problem, the paper differs
e.g.\  from standard rigged Hilbert space constructions and avoids
the introduction of nuclear spaces. The techniques are
predominantly measure theoretic and hence the Hilbert spaces
involved are separable. The results range from a Naimark type
dilation result to direct integral representations and a fairly
concrete generalized eigenvalue expansion for unbounded normal
operators.

\subsection*{Mathematics Subject Classification}
47A70 (Primary); 47A07, 47B15 (Secondary).

\subsection*{Keywords and phrases} Sesquilinear form, normal operator, generalized eigenvector, Naimark dilation, direct integral, (semi)spectral
measure.

% 47A07 Forms (bilinear, sesquilinear, multilinear)
% 47A70 (Generalized) eigenfunction expansions; rigged Hilbert spaces
% 47B15 Hermitian and normal operators (spectral measures, functional calculus, etc.)

\section{Introduction}

There is a long history of various approaches to the mathematical
clarification and justification of the well-known and
heuristically appealing formulation of Quantum Mechanics due to P.
A. M. Dirac \cite{Dirac}. Rather than try to recount this history,
we refer to the recent article \cite{GG} for references to, and a
unification of, several classical approaches based on the notion
of \emph{rigged Hilbert spaces}, and to \cite{EG:Book} for a
completely different framework of \emph{trajectory spaces}. In
Dirac's work, the notions of an ``eigenvalue'' $\lambda$ and a
corresponding ``eigenvector'' $\ket{\lambda}$ of an operator
representing a physical observable are in a key role. In general
these notions cannot be understood in their conventional
mathematical sense, and it is the task of the ancillary
mathematical theories to create a rigorous framework for their
interpretation.

The present article contributes to this line of study by providing
a relatively easily accessible setting for the analysis --
spectral in a wide sense of the word -- of a class of mathematical
objects generalizing the {\em positive operator measures} which
have been successfully used to represent physical observables in
situations where the more traditional self-adjoint operators and
their spectral measures have proved inadequate
\cite{BGL95,BLM91,Davies,Helstrom,Holevo,Ludwig}. These
generalizations, the so-called {\em sesquilinear form measures} or
{\em generalized operator measures}  were introduced in
\cite{JP:IV} in order to describe measurement situations where
only a restricted class of state preparations are available. They
also arise naturally in the quantization of classical phase
variables \cite{JP:IV}.

Despite the physical background, our study is mathematically
motivated and addresses e.g.\  some technical measurability issues
interesting in their own right.  The organization of the paper is
as follows. In Sec.~\ref{sec:defs} we introduce sesquilinear form
measures in an abstract (``test'') vector space and give an
illustrative concrete example. Sec. \ref{sec:PSFM} contains a
general representation theorem for positive sesquilinear form
measures. Some further analysis leads to a pointwise
representation theorem in Sec.~\ref{sec:pointwise}, and we relate
this result to so-called \emph{direct integral} representations in
Sec.~\ref{sec:direct}. In Sec.~\ref{sec:extension} we show that
the operations in the previous sections naturally extend from the
original test space to a larger space with a Hilbert structure.
Positive operator measures in a Hilbert space are then considered
in Sec.~\ref{sec:POM}. This section deals with a special case of
the abstract results, but we also give an alternative approach. In
Sec.~\ref{sec:normal-ops} we further specialize our results to
spectral measures of normal operators, showing that the
functionals in our general direct integral expansion indeed admit
the interpretation as generalized eigenvectors in this case. The
concluding Sec.~\ref{sec:counterexample} addresses the question of
the relationship between the spectrum of an operator and the set
of complex numbers eligible as generalized eigenvalues in some
natural sense.  An example is provided showing that even in the
case of bounded self-adjoint operators complications arise thus
indicating the need for further analysis.

When compared with much of the earlier work in this field, our
approach appears as predominantly measure theoretic and, in
particular, avoids any consideration of nuclear spaces. Inherent
in the measure theoretic setting and especially in the use of
direct integrals is the separability of the Hilbert spaces
involved. As a byproduct of the constructions of
Sec.~\ref{sec:PSFM} we obtain a generalization of the (separable
Hilbert space version of) the  Naimark dilation theorem.

\section{Preliminaries on sesquilinear form
measures}\label{sec:defs}

Let $V$ and  $W$ be vector spaces. The scalar field of all vector
spaces will be $\C$. A mapping $\Phi:V\times V\to W$ is said to be
{\em sesquilinear}, if it is antilinear (i.e., conjugate linear)
in the first and linear in the second variable. If, in addition,
$W=\C$,
 we call $\Phi$ a {\em sesquilinear form}.
 It is {\em positive}, if $V(\phi,\phi)\geq0$ for all $\phi\in V$.
 We let $S(V)$
(resp. $PS(V)$) denote the set of  sesquilinear forms (resp.
positive sesquilinear forms) on $V\times V$. It is sometimes
useful to observe that  $S(V)$ may be naturally identified with
the space of all linear maps from $V$ to $V^\times$ where
$V^\times$ is the space of antilinear functionals on $V$.
Any sesquilinear map $\Phi$ on $V\times V$ satisfies the {\em
polarization identity}
$\Phi(\phi,\psi)=\frac14\sum_{k=0}^3i^k\Phi(\phi+i^k\psi,\phi+i^k\psi).$
A positive sesquilinear form $\Phi:V\times V\to \C$ also satisfies
the equation $\Phi(\phi,\psi)=\overline{\Phi(\psi,\phi)}$ and the
\emph{Cauchy--Schwarz inequality}
$|\Phi(\phi,\psi)|^2\leq\Phi(\phi,\phi)\Phi(\psi,\psi)$ for all
$\phi,\,\psi\in V$.

Let $(\Omega,\Sigma)$ be a measurable space, i.e., $\Sigma$ is a
$\sigma$-algebra of subsets of $\Omega$.

\begin{definition}\label{def:PSFM}
Let $E:\Sigma\to S(V)$ be a mapping and denote $E(X)=E_X$ for
$X\in \Sigma$. We call $E$ a {\em sesquilinear form measure}
 if  the mapping $X\mapsto
E_X(\phi,\psi)$ is $\sigma$-additive, i.e. a complex measure, for
all $\phi,\,\psi\in V$. If, in addition, $E(\Sigma)\subset PS(V)$,
$E$ is called a {\em positive sesquilinear form measure} (or a
{\em PSF measure} or just a {\em PSFM} for short). A PSF measure
 $E:\Sigma\to PS(V)$
is called {\em strict}, if $E_\Omega(\phi,\phi)>0$ for all
$\phi\in V\setminus\{0\}$.
\end{definition}

In view of the polarization identity, in the above definition it
suffices to require that $X\mapsto E_X(\phi,\phi)$ be
$\sigma$-additive for all $\phi\in V$.  For any positive
sesquilinear form measure $E:\Sigma\to PS(V)$, the set
$N=\{\phi\in V\,|\,E_\Omega(\phi,\phi)=0\}$ is, by virtue of the
Cauchy--Schwarz inequality, a linear subspace of $V$, and the
mapping $\tilde E:\Sigma\to PS(V/N)$ which is (unambiguously)
defined by the formula
$\tilde E_X(\phi+N,\psi+N)=E_X(\phi,\psi),$
 is a strict
PSFM. It is sometimes convenient to assume that a PSFM is strict.
In view of the above observation, in many situations this
assumption does not detract essentially from generality.

Positive sesquilinear form measures may be viewed as a
generalization of some notions which we now recall. The inner
product of any Hilbert space in this paper is linear in the second
variable and usually denoted by $\<\cdot|\cdot\>$. We let
$\bddlin(H)$ denote the space of bounded linear operators on $H$,
and $\bddlin_+(H)=\{T\in\bddlin(H)\,|\,T\geq0\}$. We shall later
also encounter (the space of) the trace class (operators) on $H$,
denoted by $\bddlin^1(H)$, and its positive cone $\bddlin^1_+(H)$.
The Hilbert space of the Hilbert-Schmidt operators on $H$ is
denoted by $\bddlin^2(H)$. For a possibly unbounded operator $A$
in $H$, we denote by $\domain(A)\subset H$ its domain of
definition, and by $\range(A)=\{A\phi\,|\,\phi\in\domain(A)\}$ its
range.

\begin{definition}
Let $H$ be a Hilbert space and $E_0:\Sigma\to\bddlin_+(H)$ a
mapping. We call $E_0$ a \emph{positive operator (valued) measure}
or a \emph{POM} for short, if it weakly $\sigma$-additive, i.e.
the mapping $X\mapsto\<\phi|E_0(X)\psi\>$ is $\sigma$-additive for
all $\phi,\,\psi\in H$. If here $E_0(X)^2=E_0(X)=E_0(X)^*$ for all
$X\in \Sigma$, $E_0$ is called a projection (valued) measure. A
POM $E_0:\Sigma\to \bddlin(H)$ is called \emph{normalized} if
$E_0(\Omega)=I$, the identity operator on $H$. A normalized POM is
also called a \emph{semispectral measure} and a normalized
projection measure a \emph{spectral measure}.
\end{definition}

Every POM $E_0$ can be identified with a PSF measure $E$ by
setting $E_X(\phi,\psi):=\<\phi|E_0(X)\psi\>$.

Next we exhibit a concrete example of sesquilinear form measures.
We show that, for any weighted shift operator, there exists a
unique (shift covariant \cite{JP:V,JP:IV}) sesquilinear form
measure on the circle which is determined by a unique complex
matrix. From the structure of the matrix one can easily see when
the corresponding sesquilinear form measure is positive and
defines a POM. This is the case exactly when the shift operator is
contractive. As a byproduct we get the well-known result that the
powers of a contractive shift operator are the moment operators of
a unique semispectral measure \cite[p.\ 235]{Mlak}.

\begin{example}\label{ex:shifts}
Consider a  Hilbert space $H$ with an orthonormal basis
$(e_n)_{n\in\Z}$ and a weighted shift operator $S:\,e_n\mapsto
c_{n-1}e_{n-1}$ where $(c_n)_{n\in\Z}\subset\C$. Let
$V:=\lin(e_n)_{n\in\Z}$ and define a matrix $(c_{mn})_{m,n\in\Z}$
by $c_{mm}:=1$, $c_{mn}:=\prod_{l=m}^{n-1}c_l$ and
$c_{nm}:=\overline{c_{mn}}$ for all $m<n$. For any sesquilinear
form $\Phi:\,V\times V\to\C$ we use the formal notation
$\sum_{m,n\in\Z}\Phi(e_m,e_n)\kb{e_m}{e_n}$ which we understand as
$\sum_{m,n\in\Z}\Phi(e_m,e_n)\<\fii|e_m\>\<e_n|\psi\>=\Phi(\fii,\psi)$
for all $\fii,\,\psi\in V$.

Let $\mathcal B(\mathbb T)$ be the Borel $\sigma$-algebra of the
unit circle $\mathbb T$, and let $\nu:\,\mathcal B(\mathbb
T)\to[0,1]$ be the normalized Haar measure of $\mathbb T$. Define
a sesquilinear form measure $E_S:\,\mathcal B(\mathbb T)\to S(V)$
by
$$
E_S(X):=\sum_{m,n\in\Z}c_{mn}\int_X \lambda^{m-n}
d\nu(\lambda)\kb{e_m}{e_n},\;\;\;\;X\in \mathcal B(\mathbb T),
$$
and, for all $k\in \Z$, the $k$th moment form by
$$
E_S^{(k)}:=\sum_{m,n\in\Z}c_{mn}\int_{\mathbb T} \lambda^k
\lambda^{m-n} d\nu(\lambda)\kb{e_m}{e_n}\in S(V).
$$
Since $\int_{\mathbb T} \lambda^k \lambda^{m-n} d\nu(\lambda)=1$
when $k+m-n=0$ and $0$ otherwise, it is easy to see that, for all
$\fii,\,\psi\in V$,
$$
E_S^{(k)}(\fii,\psi)=
\begin{cases}
\<\fii|S^k\psi\>, & k>0,\\
\<\fii|(S^*)^{-k}\psi\>, & k<0,\\
\<\fii|\psi\>, & k=0,\\
\end{cases}
$$
so that the moment forms of $E_S$ can be (uniquely) extended to
the powers of $S$ and $S^*$. By analyzing the structure of $E_S$
one can derive the following result:

\begin{proposition}
$E_S$ is positive if and only if $S$ is a contraction, i.e.,
$|c_n|\le1$ for all $n\in\Z$. In this case $E_S$ has a (unique)
extension to the semispectral measure $\mathcal B(\mathbb
T)\to\mathcal L(H)$ which is a spectral measure if and only if
$|c_n|=1$ for all $n\in\Z$ (if and only if $S$ is unitary).
\end{proposition}

\begin{proof} It can be shown (see, e.g.\ \cite{JP:V} or \cite{JP:IV}) that
$E_S$ is positive (and, moreover, has a unique extension to a
semispectral measure) if and only if  the Hermitian matrix
$(c_{mn})$ is positive semidefinite. Since $(c_{mn})$ is positive
semidefinite if and only if the principal minors of $(c_{mn})$ are
nonnegative, it is easy to show that $(c_{mn})$ is positive
semidefinite if and only if $S$ is a contraction. Indeed, let
$s\in\{2,3,4,...\}$ and $\{k_1,k_2,...,k_s\}\subset\Z$ where
$k_1<k_2<...<k_s$. Since $c_{mn}=\prod_{l=m}^{n-1}c_l$, $m<n$, by
induction (see \cite[Theorem 4.1]{JP:II}), the principal minor can
be computed as
$$
\left\vert
\begin{matrix}
1 & c_{k_1,k_2} & \ldots & c_{k_1,k_s}\cr \overline{c_{k_1,k_2}} &
1 &\dots & c_{k_2,k_s}\cr \vdots & \vdots & \ddots &\vdots\cr
\overline{c_{k_1,k_s}} & \overline{c_{k_2,k_s}} &\dots & 1 \cr
\end{matrix}
\right\vert =\prod_{l=1}^{s-1}\left[1-|c_{k_l,k_{l+1}}|^2\right].
$$
Hence, $(c_{mn})$ is positive semidefinite if and only if
$|c_{mn}|=\prod_{l=m}^{n-1}|c_l|\le1$ for all $m<n$, and the claim
follows. If $E_S$ is positive, then its extension is a spectral
measure if and only if $|c_n|=1$ for all $n\in\Z$
\cite[Proposition 3]{JP:V}.
\end{proof}
\end{example}

\section{Representing positive sesquilinear form measures}\label{sec:PSFM}

For the rest of the paper, unless otherwise specified, we assume
that $(e_n)_{n=0}^{\infty}$ is a (countably infinite) Hamel basis
of $V$, indexed by $\N=\{0,\,1,\,2,\dots\}$, and $E:\Sigma\to
PS(V)$ is a positive sesquilinear form measure. We fix a sequence
of positive numbers $\alpha_n$ with
$\sum_{n=0}^\infty\alpha_n<\infty$ and write
\begin{equation}\label{eq:def-mu}
  \mu(X)=\sum_{n=0}^\infty\alpha_n
    E_X(e_n,e_n)[1+E_\Omega(e_n,e_n)]^{-1}
\end{equation}
for all $X\in \Sigma$. Then $\mu$ is a finite positive measure,
and an application of the Cauchy--Schwarz inequality shows that
for $X\in \Sigma$, $\mu(X)=0$ if and only if $E_X(\phi,\psi)=0$
whenever $\phi,\,\psi\in V$.  The Radon--Nikod\'ym theorem thus
yields for any $\phi,\,\psi\in V$ a unique element
$C(\phi,\psi)\in L^1(\mu)$ such that
$E_X(\phi,\psi)=\int_X C(\phi,\psi)\,d\mu$
for all $X\in \Sigma$. Clearly the mapping $(\phi,\psi)\mapsto
C(\phi,\psi)$ on $V\times V$ is sesquilinear, and
$C(\phi,\phi)\geq 0$ in the $L^1$-sense for all $\phi\in V$.

We let $\mathcal F$ denote the vector space of $\Sigma$-simple
$V$-valued functions on $\Omega$. If $A\subset\Omega$, $\chi_A$ is
the characteristic function of $A$, and for any
$\phi\in V$, we denote by $\phi\chi_A$ the function
$x\mapsto\chi_A(x)\phi$.

\begin{lemma} There is a unique  sesquilinear
form $\theta:\mathcal F\times\mathcal F\to \C$ satisfying
$\theta(\phi\chi_A,\psi\chi_B)=\int_{A\cap B}
C(\phi,\psi)\,d\mu$ for all $A,\,B\in\Sigma$, $\phi,\,\psi\in V$.
For any $f=\sum_{i=1}^m\phi_i\chi_{A_i}$ and
$g=\sum_{j=1}^n\psi_j\chi_{A_j}$, there holds
\begin{equation}\label{eq:theta}
  \theta(f,g)=\sum_{i=1}^m\sum_{j=1}^n\int_{A_i\cap B_j}
  C(\phi_i,\psi_j)\,d\mu,
\end{equation}
and $\theta$ is a positive sesquilinear form.
\end{lemma}

\begin{proof}
For $f,\,g\in \mathcal F$, choose representations
$f=\sum_{i=1}^m\phi_i\chi_{A_i}$ and
$g=\sum_{j=1}^n\psi_j\chi_{B_j}$ with e.g.\  all the $\phi_i$
distinct and the $A_i$ disjoint, and define $\theta(f,g)$ by the
formula (\ref{eq:theta}). Then $\theta$ is well defined, and
obvious refinement arguments yield the remaining statements.
\end{proof}

We denote $\mathcal N=\{f\in\mathcal F\,|\,\theta(f,f)=0\}$. Then
$\mathcal N$ is, by the Cauchy--Schwarz inequality, a vector
subspace of $\mathcal F$, and we get a well-defined inner product
of $\mathcal F/\mathcal N$ via
$\langle [f]|[g]\rangle=\theta(f,g)$  where e.g.\  $[f]=f+\mathcal N$.
We let $K$ denote the Hilbert space completion of this inner
product space and use the notation $\langle\cdot|\cdot\rangle$ for
the inner product of $K$.
We refer to $K$ as the \emph{associated Hilbert space} of the PSF
measure $E$ (relative to the basis $(e_n)$ and the sequence
$(\alpha_k)$).

\begin{lemma} For each $X\in \Sigma$ there is a unique bounded
linear operator $F(X):K\to K$ satisfying
$\langle [\phi\chi_A]|F(X)[\psi\chi_B]\rangle=\int_{X\cap A\cap
B}C(\phi,\psi)\,d\mu$ for all $A,\,B\in \Sigma$, $\phi,\,\psi\in
V$. Moreover, $F(X)^2=F(X)^*=F(X)$ and $F(X)[f]=[\chi_Xf]$ for all
$X\in\Sigma$, $f\in \mathcal F$.
\end{lemma}

\begin{proof} If $X\in \Sigma$, $f\in\mathcal F$
and $g\in [f]$, then
$\theta(\chi_Xf-\chi_Xg,\chi_Xf-\chi_Xg)\leq
\theta(f-g,f-g)=0,$ and so the definition $F_0(X)[f]=[\chi_Xf]$
is unambiguous. It is also clear that $\langle
F_0(X)[f]|F_0(X)[f]\rangle \leq \langle[f]|[f]\rangle$. Thus
$F_0(X)$ extends uniquely to a bounded linear map $F(X):K\to K$.
The remaining statements are immediate.
\end{proof}

\begin{lemma} The map $X\mapsto F(X)$ on $\Sigma$ is a spectral
measure.
\end{lemma}

\begin{proof} Clearly $F(\Omega)=I$, and since $\no{F(X)}\leq1$,
for weak $\sigma$-additivity it is sufficient to note that the map
$X\mapsto\langle [\phi\chi_A]|F(X)[\psi\chi_B]\rangle=\int_{X\cap
A\cap B} C(\phi,\psi)\,d\mu$  on $\Sigma$ is $\sigma$-additive
for all for all $A,\,B\in \Sigma$, $\phi,\,\psi\in V$.
\end{proof}

We now define $J:V\to K$ by the formula $J\phi=[\phi\chi_\Omega]$.
Then $J$ is a linear map. We collect and complement the above
arguments in the following theorem.

\begin{theorem}\label{th:dilation} Let $E:\Sigma\to PS(V)$ be a PSFM.

{\rm(a)} There is a Hilbert space $K$ with a spectral measure
$F:\Sigma\to\bddlin(K)$ and a linear map $J:V\to K$ such that
$\langle
J\phi|F(X)J\psi\rangle= E_X(\phi,\psi)$ for all $X\in\Sigma$ and
$\phi,\,\psi\in V$, and moreover, the linear span of the set
$\{F(X)J\phi\,|\,X\in\Sigma,\,\phi\in V\}$ is dense in $K$.

{\rm(b)} This representation of $E$ is essentially unique in the
sense that if the triple $(K_1,F_1,J_1)$ gives another
representation with these properties, there is a unique unitary
map $U:K\to K_1$ such that $UF(X)J\phi=F_1(X)J_1\phi$ for all
$X\in\Sigma$, $\phi\in V$; in particular, $UJ\phi=J_1\phi$ for all
$\phi\in V$. Moreover, $UF(X)=F_1(X)U$ for all $X\in \Sigma$.

{\rm(c)} In the situation of (a), $J$ is injective if and only if
$E$ is strict.
\end{theorem}

\begin{proof}
(a) In the above construction,
$\langle
J\phi|F(X)J\psi\rangle= \langle
[\phi\chi_\Omega]|F(X)[\phi\chi_\Omega]\rangle=\int_XC(\phi,\psi)\,d\mu
=E_X(\phi,\psi).$ The density statement is also an immediate
consequence of the construction.

(b) The uniqueness of $U$ is clear since it is determined on a
dense subspace of $K$. On the other hand, if
$X_1,\dots,X_n\in\Sigma$ and $\phi_1,\dots,\phi_n\in V$, then
\[\begin{split}
&\no{\sum_{i=1}^nF(X_i)J\phi_i}^2=\sum_{i=1}^n\sum_{j=1}^n\<F(X_i)J\phi_i|F(X_j)J\phi_j\>\\
&=\sum_{i=1}^n\sum_{j=1}^n\<J\phi_i|F(X_i\cap X_j)J\phi_j\>=\sum_{i=1}^n\sum_{j=1}^n E_{X_i\cap
X_j}(\phi_i,\phi_j)=\no{\sum_{i=1}^nF_1(X_i)J_1\phi_i}^2.
\end{split}\]
Thus there is a well-defined isometric linear map sending each
$\sum_{i=1}^nF(X_i)J\phi_i$ to $\sum_{i=1}^nF_1(X_i)J_1\phi_i$,
and this map extends by continuity to a unitary $U:K\to K_1$. In
particular, $UJ\phi=UF(\Omega)J\phi=F_1(\Omega)J_1\phi=J_1\phi$
for all $\phi\in V$. Moreover, $UF(X)F(Y)J\phi=UF(X\cap Y)J\phi=
F_1(X\cap Y)J_1\phi=F_1(X)F_1(Y)J_1\phi=F_1(X)UF(Y)J\phi$ for all
$X,\,Y\in \Sigma$, $\phi\in V$, from which the equation
$UF(X)=F_1(X)U$ follows.

(c) Suppose first that the triple $(K,F,J)$ is obtained by the
measure theoretic construction preceding the theorem. Then
\begin{equation}\label{ekvatsione}
\no{J\phi}^2=\no{[\phi\chi_\Omega]}^2
=\theta(\phi\chi_\Omega,\phi\chi_\Omega)=\int_\Omega
C(\phi,\phi)\,d\mu=E_\Omega(\phi,\phi),
\end{equation}
so in particular $J\phi$ vanishes if and only if
$E_{\Omega}(\phi,\phi)$ does, and the claim follows. In the case
of another triple $(K_1,F_1,J_1)$, let $U:K\to K_1$ be as in (b).
Since $UJ\phi=J_1\phi$ and $U$ is bijective, $J\phi=0$ if and only
if $J_1\phi=0$.
\end{proof}

\begin{remark}\label{rem:measure independence}
The uniqueness part of the above result shows, in particular, that
the choice of the basis $(e_n)$ and the sequence $(\alpha_n)$ does
not essentially influence the resulting structure. In fact $\mu$
could be replaced by any finite positive measure $\nu$ such that
every complex measure $X\mapsto E_X(\phi,\psi)$ is absolutely
continuous with respect to $\nu$.
\end{remark}

\begin{remark}\label{rem:Stinespring}
Let $H$ be a separable Hilbert space, $(e_n)_{n=0}^{\infty}$
 an orthonormal basis of $H$, and $V:=\lin(e_n)_{n=0}^{\infty}$
 its linear span. Suppose that
$E_0:\Sigma\to\bddlin(H)$ is a semispectral measure and let
$E:\Sigma\to PS(V)$ be the PSFM defined by $E_X(\phi,\psi)
=\langle\phi|E_0(X)\psi\rangle$ for $X\in \Sigma$, $\phi,\,\psi\in
V$. Retaining the notation of the general case, now $J:V\to K$ is
isometric, because $\no{J\phi}^2=E_\Omega(\phi,\phi)=\no{\phi}^2$
by (\ref{ekvatsione}).

In this case the spectral measure $F:\Sigma\to \bddlin(K)$ is the
minimal Naimark dilation of $E_0$ (see e.g.\  \cite{Paulsen}). It
follows form an observation made in \cite[p.\ 171]{LY04} that the
original semispectral measure $E_0$ is a spectral measure if and
only if $J(V)$ is dense in $K$, or equivalently,
$(Je_n)_{n=0}^{\infty}$ is an orthonormal basis of $K$.
\end{remark}
\section{Pointwise representation}\label{sec:pointwise}

In the previous section we obtained the representation
$
  E_X(\phi,\psi)=\int_X C(\phi,\psi)d\mu,
$
where $C:V\times V\to L^1(\mu)$ is sesquilinear and
$C(\phi,\phi)\geq 0$ in the $L^1(\mu)$ sense. As a step towards
the main result of this section, let us now provide a pointwise
version of this formula. We use the notion of $\mu$-measurability
as in \cite{DS} and often call it just measurability.
 Since
$\mu$ is a finite measure, for scalar function this simply means
measurability with respect to the Lebesgue extension of $\Sigma$
with respect to $\mu$.
 We say that
$\Omega\owns\omega\mapsto C_{\omega}\in PS(V)$ is a
($\mu$-)\emph{measurable family of forms} if $\omega\mapsto
C_{\omega}(\phi,\psi)$ is $\mu$-measurable for all $\phi,\,\psi\in
V$.

\begin{lemma}\label{lem:pointwise}
For a PSF measure $E:\Sigma\to PS(V)$, there is a measurable
family of forms $\Omega\owns\omega\mapsto C_{\omega}\in PS(V)$,
such that for all $\phi,\,\psi\in V$, the function $\omega\mapsto
C_{\omega}(\phi,\psi)$ is a representative of $C(\phi,\psi)\in
L^1(\mu)$.
\end{lemma}

\begin{proof}
For all $m,\,n\in\N$, let us pick a function representative
$g_{mn}\in C(e_m,e_n)$. For every $\omega\in\Omega$ and any
$\phi,\,\psi\in V$, we define
\begin{equation}\label{eq:pointwise-form}
  C_{\omega}(\phi,\psi):=\sum_{m,n=0}^{\infty}
    \bar{a}_m g_{mn}(\omega)b_n,
\end{equation}
where $\phi=\sum_{m=0}^{\infty}a_m e_m$,
$\psi=\sum_{n=0}^{\infty}b_n e_n$ are the unique expansions in the
Hamel basis $(e_n)_{n=0}^{\infty}$, and only finitely many of the
coefficients $a_m,\,b_n\in\C$ are non-zero. Then it is immediate
that $C_{\omega}:V\times V\to\C$ is a sesquilinear form, and the
(measurable) function $\omega\mapsto C_{\omega}(\phi,\psi)$ is a
representative of $C(\phi,\psi)$. In particular
$C_{\omega}(\phi,\phi)\geq 0$ for $\mu$-a.e.\ $\omega\in\Omega$.
If we only consider the countable set
$W:=\lin_{\Q+i\Q}(e_n)_{n=0}^{\infty}$, we can choose a single
$\mu$-null set $Z\subset\Omega$ such that
$C_{\omega}(\phi,\phi)\geq 0$ for all $\omega\in\Omega\setminus Z$
and $\phi\in W$. But for a general $\phi\in V$, we may approximate
the finitely many non-zero coefficient $a_m$ in
(\ref{eq:pointwise-form}) by complex rationals, getting the same
result for all $\phi\in V$. Thus $C_{\omega}\in PS(V)$ for all
$\omega\in\Omega\setminus Z$. If we redefine
$C_{\omega}(\phi,\psi):=0$ for $\omega\in Z$ (which is achieved by
changing, if necessary, the functions $g_{mn}$ to have the value
zero in $Z$), we have $C_{\omega}\in PS(V)$ for all
$\omega\in\Omega$, and $\omega\mapsto C_{\omega}(\phi,\psi)$ is
still a representative of $C(\phi,\psi)$ for all $\phi,\,\psi\in
V$.
\end{proof}
We now introduce some notational conventions which depend on the
choice of the fixed Hamel basis $(e_n)_{n=0}^\infty$ of $V$. If
$\phi,\,\psi\in V$ have the basis expansions
$\phi=\sum_{n=0}^\infty a_ne_n$ and $\psi=\sum_{n=0}^\infty
b_ne_n$ (with only finitely many non-zero terms), we write
$\<\phi|\psi\>=\sum_{n=0}^\infty\bar{a}_nb_n.$
Then $V$ becomes an inner-product space, and each $\phi\in V$
gives rise to the linear form $\psi\mapsto\<\phi|\psi\>$, denoted
by $\bra{\phi}$, and to the antilinear form
$\psi\mapsto\<\psi|\phi\>$, denoted by $\ket{\phi}$. Clearly, the
space of all linear functionals on $V$, i.e., the algebraic dual
$V^*$ of $V$, is in a bijective antilinear correspondence with the
vector space of all complex sequences $(d_n)$, when $d=(d_n)$ is
made to correspond to the functional $ \psi=\sum_{n=0}^\infty
b_ne_n\mapsto \sum_{n=0}^\infty\bar{d}_nb_n$, denoted also by
$\bra{d}$ in the sequel; we may write this as
$\<d|\psi\>=\sum_{n=0}^\infty\bar{d}_nb_n.$
In this kind of situations we also allow the notation
$\<\psi|d\>=\sum_{n=0}^\infty\bar{b}_nd_n,$
so that  $\<\psi|d\>=\overline{\<d|\psi\>}$. These notations are
consistent when we use the identification of a vector
$\phi=\sum_{n=0}^\infty a_ne_n\in V$ with the sequence $(a_n)$ of
its coefficients. Extending this identification, we  sometimes
consider $V$ embedded (linearly) in the space of all sequences
$(a_n)_{n=0}^\infty$ and identify a sequence $(a_n)_{n=0}^\infty$
with the formal series $\sum_{n=0}^\infty a_ne_n$. In particular
note that $\bra{e_n}$ is the linear functional on $V$ which maps
$\phi=\sum_{k=0}^{\infty}a_k e_k$ into $a_n$. It is sometimes
convenient to denote by the formal series $\sum_{n=0}^\infty
a_n\bra{e_n}$ the element of $V^*$ corresponding in our convention
to the  sequence $(\bar{a}_n)$. If $\phi=\sum_{n=0}^\infty
a_ne_n\in V$, i.e., the series is not just formal, then
$\bra{\phi}=\sum_{n=0}^\infty\bar{a}_n\bra{e_n}.$

If $\bra{d_1},\,\bra{d_2}\in V^*$, we denote by
 $\ket{d_1}\bra{d_2}$ the sesquilinear form
  $(\phi,\psi)\mapsto\<\phi|d_1\>\<d_2|\psi\>$; it can
equivalently be viewed as the antilinear map
$\phi\mapsto\<\phi|d_1\>\bra{d_2}$ from $V$ to $V^*$. Note that
$\ket{d}\bra{d}\in PS(V)$ for any $\bra{d}\in V^*$.

In the sequel, we make various measurability assertions concerning
vector-valued functions. When saying that a $V$-valued function
$\omega\mapsto F(\omega)$ is measurable, unless otherwise
specified, we mean \emph{weak measurability} with respect to the
dual pair $\<V,V^*\>$, i.e., that the scalar-valued functions
$\omega\mapsto\<d|F(\omega)\>$ are measurable for all $\bra{d}\in
V^*$.

\begin{lemma}\label{lem:measurability}
A $V$-valued function $\omega\mapsto F(\omega)$ is measurable if
and only if all the coordinate functions
$\omega\mapsto\<e_n|F(\omega)\>$, $n\in\N$, are measurable.
\end{lemma}

\begin{proof}
One direction is clear, since $\bra{e_n}\in V^*$. But if the
coordinate functions are measurable and
$\bra{d}=\sum_{k=0}^{\infty}c_k\bra{e_k}\in V^*$, then also
$\<d|F(\omega)\>=\sum_{k=0}^{\infty}c_k\<e_k|F(\omega)\>$ is
measurable as the sum of an everywhere convergent series of
measurable functions.
\end{proof}

\begin{lemma}\label{lem:diag}
For a $\mu$-measurable family of forms $C_{\omega}\in PS(V)$,
$\omega\in\Omega$, there exist $\mu$-measurable mappings
$\omega\mapsto n(\omega)\in\N\cup\{\infty\}$ and $\omega\mapsto
g_k(\omega)\in V$, $k\in\N$, such that for all $\omega\in\Omega$:
\begin{itemize}
  \item $C_{\omega}(g_k(\omega),g_{\ell}(\omega))=
    \delta_{k\ell}\chi_{\{\omega'|k<n(\omega')\}}(\omega)$.
  \item $\displaystyle C_{\omega}(\phi,\psi)
  =\sum_{k=0}^{n(\omega)-1}C_{\omega}(\phi,g_k(\omega))
   C_{\omega}(g_k(\omega),\psi)$ for all $\phi,\,\psi\in V$, and
   only finitely many terms are non-zero for fixed $\phi$ and
   $\psi$ even when $n(\omega)=\infty$.
   \item $\omega\mapsto C_{\omega}(g_k(\omega),\phi)$ is
   $\mu$-measurable for every $\phi\in V$.
\end{itemize}
\end{lemma}

\begin{proof}
For every $\omega\in\Omega$, we denote
$\mathcal{N}_{\omega}:=\{\phi\in V\,|\,C_{\omega}(\phi,\phi)=0\}$.
In complete analogy with the space $\mathcal{N}$ of the previous
section, we see that this is a linear subspace of $V$, and on
$V/\mathcal{N}_{\omega}$ we get a well-defined inner-product
$\langle [\phi]_{\omega}|[\psi]_{\omega}\rangle_{\omega}
:=C_{\omega}(\phi,\psi)$, where e.g.\ \
$[\phi]_{\omega}:=\phi+\mathcal{N}_{\omega}$.

We now apply the Gram--Schmidt procedure on
$(V/\mathcal{N}_{\omega},\langle\cdot|\cdot\rangle_{\omega})$,
starting from the spanning sequence
$([e_n]_{\omega})_{n=0}^{\infty}$. Instead of doing this at the
completely abstract level, however, we work with concrete
representatives in $V$ in order to keep track of the
$\mu$-measurability of our operations. We also do not discard
possible zero-vectors in the first place.

Denote $\{\phi\}_{\omega}^0:=C_{\omega}(\phi,\phi)^{-1/2}\phi$ if
$[\phi]_{\omega}\neq 0$ and $\{\phi\}_{\omega}^0:=0$ otherwise.
Then let
\[
  f_0(\omega):=\{e_0\}_{\omega}^0,\qquad
  f_n(\omega):=\big\{e_n-\sum_{k=0}^{n-1}
    C_{\omega}(f_k(\omega),e_n)f_k(\omega)\big\}_{\omega}^0,
  \quad n=1,2,\ldots
\]
It is easily seen that
$f_n(\omega)\in\lin\{e_0,\ldots,e_n\}$ and
$
  \lin\{f_0(\omega),\ldots,f_n(\omega)\}+\mathcal{N}_{\omega}
  =\lin\{e_0,\ldots,e_n\}+\mathcal{N}_{\omega}
$
for all $n\in\N$ and $\omega\in\Omega$.
Moreover, the functions $\omega\mapsto f_n(\omega)$ are
$\mu$-measurable.

Next, we remove the possible zero-vectors in a measurable way: Let
$n_0(\omega):=-1$ and
\[\begin{split}
  n_k(\omega)&:=\inf\{n\in\N\,|\,n>n_{k-1},f_n(\omega)\neq 0\}
   \in\N\cup\{\infty\},\qquad k\in\N,\\
  n(\omega)&:=1+\sup\{k\in\N\cup\{-1\}\,|\,n_k(\omega)<\infty\}
   \in\N\cup\{\infty\},
\end{split}\]
where $\inf\emptyset:=\infty$. Finally, we set
$g_k(\omega):=f_{n_k(\omega)}(\omega)$ for $0\leq k<n(\omega)$.
These are the non-zero vectors from the Gram--Schmidt procedure,
and hence $([g_k(\omega)]_{\omega})_{k=0}^{n(\omega)-1}$ is an
orthonormal Hamel basis of
$(V/\mathcal{N}_{\omega},\langle\cdot|\cdot\rangle_{\omega})$.
This implies the first two claims after setting $g_k(\omega):=0$
for $k\geq n(\omega)$. The last claim follows from the formula
\[
  C_{\omega}(g_k(\omega),\phi)
  =\sum_{j=0}^{n_k(\omega)-1}\bar{f}_{n_k(\omega),j}(\omega)
   C_{\omega}(e_j,\phi),
\]
where $f_{n_k(\omega),j}(\omega)$ stands for the $j$th coordinate
of $f_{n_k(\omega)}(\omega)$, and from the measurability of sums
and products of measurable functions.
\end{proof}

In the above representation, $\phi\mapsto
C_{\omega}(g_k(\omega),\phi)$ is a linear functional on $V$. There
is a unique complex sequence $d_k(\omega)$ such that this
functional equals $\<d_k(\omega)|\in V^*$. By the measurability of
$V^*$-valued functions, we understand the weak$^*$ measurability,
i.e., the measurability of their pointwise duality pairings with
all $\phi\in V$. By an argument similar to, but even easier than,
Lemma~\ref{lem:measurability}, this is seen to be equivalent to
the measurability of all the coordinate functions, i.e., it
suffices to test the pairings with all $\phi=e_n$, $n\in\N$. For
$\omega\mapsto \bra{d_k(\omega)}$ this measurability condition is
precisely the last claim of the previous lemma. With these
remarks, and a combination of Lemmas~\ref{lem:pointwise} and
\ref{lem:diag}, we have the following result, in which
$E:\Sigma\to PS(V)$ is a PSFM, and all our earlier choices and
notations apply.

\begin{theorem}\label{th:diag}
There are $\mu$-measurable mappings $\omega\mapsto
n(\omega)\in\N\cup\{\infty\}$, $\omega\mapsto g_k(\omega)\in V$
and $\omega\mapsto \bra{d_k(\omega)}\in V^*$, $k\in\N$, such that
$\langle d_k(\omega)|g_{\ell}(\omega)\rangle=
\delta_{k\ell}\chi_{\{\omega'|k<n(\omega')\}}(\omega)$ for all
$\omega\in\Omega$, and the following representation holds for all
$\phi,\,\psi\in V$:
\[
  E_X(\phi,\psi)=\int_X\sum_{k=0}^{n(\omega)-1}
   \langle\phi|d_k(\omega)\rangle
   \langle d_k(\omega)|\psi\rangle d\mu(\omega),
\]
where only finitely many terms in the sum are non-zero for each
$\omega$, even when $n(\omega)=\infty$. In particular,
\begin{equation}\label{eq:isometry}
  E_X(\phi,\phi)=\int_X\sum_{k=0}^{n(\omega)-1}
   |\langle d_k(\omega)|\phi\rangle|^2 d\mu(\omega).
\end{equation}
\end{theorem}

\section{Relation to direct integrals}\label{sec:direct}

In this section we make some comments on the relation of
Theorem~\ref{th:diag} to the direct integral representations which
are popular in some of the related literature. We use the
following notational conventions. As usual, $\ell^2$ is the
Hilbert space of square summable sequences $(a_n)_{n=0}^\infty$,
and $L^2(\Omega,\mu,\ell^2)$ is the Hilbert space consisting of
the ($\mu$-equivalence classes of) $\mu$-measurable
$\ell^2$-valued functions $f$ for which the function
$\omega\mapsto\no{f(\omega)}^2$ is $\mu$-integrable. For any
$k\in\N\cup\{\infty\}$ we denote
$\ell^2_k=\{(a_n)\in\ell^2\,|\,a_n=0 \text{ for all } n\geq
k\}.$ In particular, $\ell^2_0=\{0\}$ and $\ell^2_\infty=\ell^2$.

Let now $\omega\mapsto n(\omega)$ be a $\mu$-measurable map from
$\Omega$ into $\N\cup\{\infty\}$. We denote by
$L^2_{n(\cdot)}(\Omega,\mu,\ell^2)$ the subset of
$L^2(\Omega,\mu,\ell^2)$ consisting of those $f$ for which
$f(\omega)\in \ell^2_{n(\omega)}$ $\mu$-almost everywhere. A
routine argument based on the fact that a sequence converging in
$L^2(\Omega,\mu,\ell^2)$, and hence in $\mu$-measure, contains an
almost everywhere convergent subsequence shows that
$L^2_{n(\cdot)}(\Omega,\mu,\ell^2)$ is a closed subspace of
$L^2(\Omega,\mu,\ell^2)$.

The above construction of the space
$L^2_{n(\cdot)}(\Omega,\mu,\ell^2)$ is a relatively simple way of
making rigorous the heuristic idea of an $L^2$ space of functions
taking their pointwise values in Hilbert spaces of different
dimensions. This type of space is often referred to as the
\emph{direct integral} of Hilbert spaces and denoted by
$
  \int_{\Omega}^{\oplus}\ell^2_{n(\omega)}d\mu(\omega).
$

We now return to Theorem \ref{th:diag} and its notation; in
particular, let $\omega\mapsto n(\omega)$ henceforth stand for the
fixed map appearing in that Theorem. If $\phi\in V$, define
$J_1\phi$ to be the element
$\omega\mapsto(\<d_k(\omega)|\phi\>)_{k=0}^{\infty}$ of
 $L^2_{n(\cdot)}(\Omega,\mu,\ell^2)$. This produces a linear map $J_1:V\to L^2_{n(\cdot)}(\Omega,\mu,\ell^2)$. For $X\in \Sigma$,
 let $F_1(X)$ be the restriction of the multiplication map
 $f\mapsto \chi_Xf$ to the invariant
 subspace $L^2_{n(\cdot)}(\Omega,\mu,\ell^2)$ of $L^2(\Omega,\mu,\ell^2)$.
 Then $F_1:\Sigma\to \bddlin(L^2_{n(\cdot)}(\Omega,\mu,\ell^2))$ is a spectral measure.

\begin{theorem}\label{th:dirint}
The linear span of the set
$\{F_1(X)J_1\phi\,|\,X\in\Sigma,\,\phi\in V\}$ is dense in
$K_1:=L^2_{n(\cdot)}(\Omega,\mu,\ell^2)$, and
$\<J_1\phi|F_1(X)J_1\psi\>=E_X(\phi,\psi)$
for any $X\in \Sigma$ and $\phi,\,\psi\in V$. Thus the triple
$(K_1,F_1,J_1)$ is unitarily equivalent to the triple $(K,F,J)$ in
the sense of Theorem~\textup{\ref{th:dilation}}.
\end{theorem}

\begin{proof}
To prove the density statement, let
$G=[\omega\mapsto(G_k(\omega))_{k=0}^{\infty}]\in
L^2_{n(\cdot)}(\Omega,\mu,\ell^2)$ be  orthogonal to
$F_1(X)J_1e_m$ for all $m\in\N$ and $X\in\Sigma$. This means that
\[
  \int_X\sum_{k=0}^{n(\omega)-1}\bar{G}_k(\omega)
    \langle d_k(\omega)|e_m\rangle d\mu(\omega)=0\qquad
  \text{for all }X\in\Sigma,\, m\in\N.
\]
But this implies that
$
  \sum_{k=0}^{n(\omega)-1}\bar{G}_k(\omega)
    \langle d_k(\omega)|e_m\rangle =0
$
for all $m\in\N$ 
and $\mu$-a.e.\ $\omega\in\Omega,$
where the exceptional $\mu$-null set, say $Z$, may be taken
independent of $m\in\N$.
Recall that $g_{\ell}(\omega)\in V=\lin(e_m)_{m=0}^{\infty}$. By
linearity, we obtain
$
  \sum_{k=0}^{n(\omega)-1}\bar{G}_k(\omega)
    \langle d_k(\omega)|g_{\ell}(\omega)\rangle=0
    $
for all $\ell\in\N$ and all 
$ \omega\in\Omega\setminus Z.$
But the left-hand side is $\bar{G}_{\ell}(\omega)$, and hence we
have shown that $G$ vanishes as an element of
$L^2_{n(\cdot)}(\Omega,\mu,\ell^2)$. This proves our first claim.
If $X\in \Sigma$ and $\phi,\,\psi\in V$, then
$$\<J_1\phi|F_1(X)J_1\psi\>=\int_X\sum_{k=0}^{n(\omega)-1}
   \<\phi|d_k(\omega)\>\<d_k(\omega)|\psi\>
   d\mu(\omega)=E_X(\phi,\psi)$$
by Theorem~\ref{th:diag}. The asserted unitary equivalence now
follows from Theorem~\ref{th:dilation}.
\end{proof}

As a consequence of Theorem~\ref{th:dirint} and
Remark~\ref{rem:Stinespring}, we have:

\begin{corollary}\label{cor:Dirac-basis}
Let $E_X(\phi,\psi)=\<\phi|E_0(X)\psi\>$ for a semispectral
measure $E_0$. Then $E_0$ is a spectral measure if and only if
\begin{equation}\label{eq:Dirac-basis}
  \omega\mapsto(\<d_k(\omega)|e_n\>)_{k=0}^{\infty},\
  n\in\N,\text{ is an orthonormal basis of }
  L^2_{n(\cdot)}(\Omega,\mu,\ell^2).
\end{equation}
\end{corollary}

\section{Extension of the test vector space}\label{sec:extension}

Our considerations so far have taken place in the space $V$, which
only has an algebraic vector space structure. We now show that
 the operations on $V$ that we have been dealing
with actually extend to a Hilbert space completion of this initial
space that we started with. Let us denote
\[
  H:=\Big\{\psi=\sum_{n=0}^{\infty}b_ne_n\,\Big|\,
    \|\psi\|_H^2:=\sum_{n=0}^{\infty}|b_n|^2<\infty\Big\}.
\]
Then $H$ is a Hilbert space with the orthonormal basis
$(e_n)_{n=0}^{\infty}$, whose linear span is $V$.

Recall the definition of the measure $\mu$ from (\ref{eq:def-mu}).
We use the same quantities from this definition to introduce the
linear operator $\Lambda:H\to H$, given by
\begin{equation}\label{eq:Lambda}
  \Lambda:=\sum_{n=0}^{\infty}\beta_n
    |e_n\>\<e_n|,\qquad
  \beta_n:=\frac{\alpha_n}{1+E_{\Omega}(e_n,e_n)}.
\end{equation}
This operator is positive and injective, and its range is dense in
$H$. We shall also use the functional calculus of $\Lambda$, i.e.,
the operators
$\eta(\Lambda):=\sum_{n=0}^{\infty}\eta(\beta_n)|e_n\>\<e_n|$ for
functions $\eta:\R_+\to\C$. Note that $\eta(\Lambda)V\subset V$.

Let us now consider the form
$C_{\omega}(\Lambda^{1/2}\cdot,\Lambda^{1/2}\cdot)\in PS(V)$. We
define its \emph{trace} by
\[
  \tr(C_{\omega}(\Lambda^{1/2}\cdot,\Lambda^{1/2}\cdot))
  :=\sum_{n=0}^{\infty}C_{\omega}(\Lambda^{1/2}e_n,\Lambda^{1/2}e_n)
  =\sum_{n=0}^{\infty}\beta_n C_{\omega}(e_n,e_n).
\]
Integration over an arbitrary $X\in\Sigma$ gives
\[
  \int_X\tr(C_{\omega}(\Lambda^{1/2}\cdot,\Lambda^{1/2}\cdot))
    d\mu(\omega)
  =\sum_{n=0}^{\infty}\beta_n E_X(e_n,e_n)
  =\mu(X).
\]
This means that
$\tr(C_{\omega}(\Lambda^{1/2}\cdot,\Lambda^{1/2}\cdot))=1$ for
$\mu$-a.e.\ $\omega\in\Omega$, say for $\omega\in\Omega\setminus
Z$, where $\mu(Z)=0$.

Let us momentarily restrict ourselves to a finite dimensional
space $V_N:=\lin(e_n)_{n=0}^N$, which we make into a Hilbert space
such that $(e_n)_{n=0}^N$ is an orthonormal basis. The restriction
of $C_{\omega}(\Lambda^{1/2}\cdot,\Lambda^{1/2}\cdot)$ to
$V_N\times V_N$ belongs to $PS(V_N)$ and defines a positive
operator, say $T_N(\omega)$, on $V_N$ in a natural way. The
computation of the previous paragraph implies that the (usual)
trace of $T_N(\omega)$, and hence its norm in $\bddlin(V_N)$, is
at most $1$, for all $\omega\in\Omega\setminus Z$. This uniform
estimate, and the density of $V=\bigcup_{n=0}^{\infty}V_N$ in $H$,
imply by a standard limiting argument the existence of an operator
$T(\omega)\in\bddlin(H)$ such that
$\<\phi|T(\omega)\psi\>=C_{\omega}(\Lambda^{1/2}\phi,\Lambda^{1/2}\psi)$
for all $\phi,\,\psi\in V$, and this operator is positive with
trace $1$. For definiteness, let us define
$T(\omega):=0\in\bddlin(H)$ for all $\omega$ in the exceptional
set $Z$.

We now have the equations
\[
  C_{\omega}(\phi,\psi)
  =\<\Lambda^{-1/2}\phi|T(\omega)\Lambda^{-1/2}\psi\>,
  \qquad
  E_X(\phi,\psi)=\int_X
   \<\Lambda^{-1/2}\phi|T(\omega)\Lambda^{-1/2}\psi\>d\mu(\omega)
\]
for all $\phi,\,\psi\in V$, but we can see that the right-hand
sides are meaningful for all
$\phi,\,\psi\in\domain(\Lambda^{-1/2})=\range(\Lambda^{1/2})=:H_1$,
where we define
\begin{equation}\label{eq:Hgamma}
  H_{\gamma} :=\Big\{\phi=\sum_{n=0}^{\infty}a_n e_n\,\Big|\,
  \|\phi\|_{H_\gamma}^2:=
  \sum_{n=0}^{\infty}\frac{|a_n|^2}{\beta_n^{\ \gamma}}<\infty\Big\},
  \qquad\gamma\in\R.
\end{equation}

Since $\sum_{k=0}^{n(\omega)-1}|\<d_k(\omega)|\phi\>|^2
=C_{\omega}(\phi,\phi)=\<\Lambda^{-1/2}\phi|T(\omega)\Lambda^{-1/2}\phi\>\leq\|\Lambda^{-1/2}\phi\|^2$,
we also find that the functionals $\bra{d_k(\omega)}\in V^*$ in
fact extend by continuity to continuous linear functionals on
$H_1$. Observe that under the natural duality of sequence spaces
we have $H_1^*\eqsim H_{-1}$. In particular, the coordinate
sequences satisfy
\begin{equation}\label{eq:sizeOfDk}
  (\<d_k(\omega)|e_n\>\beta_n^{1/2})_{n=0}^{\infty}\in\ell^2.
\end{equation}

We collect the results from the above discussion in the following:

\begin{proposition}
The sesquilinear forms $E_X,C_{\omega}\in PS(V)$, $X\in\Sigma$,
$\omega\in\Omega$, and the functionals $\bra{d_k(\omega)}\in V^*$
extend continuously to the Hilbert space completion $H_1$ of $V$
\textup{(}defined in \textup{(\ref{eq:Hgamma}))}, which is
associated with the PSFM $E$.

On the larger Hilbert space $H=H_0$, which is canonically related
to $V$ and independent of $E$, they admit representations as
unbounded forms with common dense domain:
\[
  E_X(\phi,\psi)=\int_X C_{\omega}(\phi,\psi)d\mu(\omega),\quad
  C_{\omega}(\phi,\psi)
  =\<\Lambda^{-1/2}\phi|T(\omega)\Lambda^{-1/2}\psi\>,\quad
  \phi,\,\psi\in\range(\Lambda^{1/2})=H_1,
\]
with an injective $\Lambda\in\bddlin^1_+(H)$ \textup{(}defined in
\textup{(\ref{eq:Lambda}))} and
\[
  T(\omega)=\sum_{k=0}^{n(\omega)-1}
  \Lambda^{1/2}\ket{d_k(\omega)}\bra{d_k(\omega)}\Lambda^{1/2}
  =:\sum_{k=0}^{n(\omega)-1}\ket{h_k(\omega)}\bra{h_k(\omega)}
  \in\bddlin^1_+(H),
\]
where $h_k(\omega)\in H$ is defined with the help of the Riesz
representation theorem in terms of the linear functional
$\<h_k(\omega)|\phi\>:=\<d_k(\omega)|\Lambda^{1/2}\phi\>$,
$\phi\in H$. Moreover,
\[
  \tr\,T(\omega)=\sum_{k=0}^{n(\omega)-1}
    \|h_k(\omega)\|^2_H=1,\quad
  \rank T(\omega)=n(\omega)
  \qquad\text{for a.e. }\omega\in\Omega,
\]
and the maps $\Omega\ni\omega\mapsto T_{\omega}\in\bddlin^1(H)$
and $\Omega\ni\omega\mapsto h_k(\omega)\in H$ are B\^ochner
$\mu$-measurable.
\end{proposition}

\begin{proof}[Rest of the proof]
What remains to show is in the last two lines of the assertions.
Concerning the rank of $T(\omega)$, we know that
$
  \delta_{k\ell}\chi_{\{\omega'\in\Omega|k\leq
  n(\omega')\}}(\omega)
  =\<d_k(\omega)|g_{\ell}(\omega)\>
  =\<h_k(\omega)|\Lambda^{-1/2}g_{\ell}(\omega)\>,
$
where $g_{\ell}(\omega)\in V\subset\range(\Lambda^{1/2})$ so that
$\Lambda^{-1/2}g_{\ell}(\omega)\in H$. This shows that the vectors
$h_k(\omega)$, $0\leq k<n(\omega)$, are linearly independent.

As for measurability, we know that
$\omega\mapsto\<h_k(\omega)|\psi\>=\<d_k(\omega)|\Lambda^{-1/2}\psi\>$
is measurable for all $\psi\in V$. By the density of $V$ in $H$,
the measurability extends to all $\psi\in H$, and the weak
measurability thus established is equivalent to the (strong)
B\^ochner measurability in the separable Banach space $H$ by the
Pettis measurability theorem (see e.g.~\cite[Theorem~II.1.2]{DU}
or \cite[Theorem~III.6.11]{DS}). This implies the measurability of
the finite sums
$\sum_{k=0}^{\min(N,n(\omega)-1)}\ket{h_k(\omega)}\bra{h_k(\omega)}$
convergent to $T(\omega)$ (pointwise in the norm of
$\bddlin^1(H)$), which is hence measurable also.
\end{proof}

\section{Positive operator measures}\label{sec:POM}

In this section we consider an important special case of the above
theory, where the PSFM is defined on the whole Hilbert space $H$
from the beginning. Now $(e_n)_{n=0}^{\infty}$ is an orthonormal
basis of $H$, but we continue to denote
$V=\lin(e_n)_{n=0}^{\infty}$ as in the previous sections.
  Every POM $E_0:\Sigma\to \bddlin(H)$ can be
identified with a PSF measure $E$ by setting
$E_X(\phi,\psi):=\<\phi|E_0(X)\phi\>$. Thus the general results
give:

\begin{proposition}\label{prop:pomrep}
Given a POM $E_0:\Sigma\to\bddlin(H)$ and defining a measure $\mu$
by \textup{(\ref{eq:def-mu})}, there exists the following
representation:
\[
  \<\phi|E_0(X)\psi\>
  =\int_X\<\Lambda^{-1/2}\phi|T(\omega)
  \Lambda^{-1/2}\psi\>d\mu(\omega),\qquad
  \phi,\,\psi\in\range(\Lambda^{1/2}),
\]
where $\Lambda\in\bddlin^1_+(H)$ is injective,
\[
T(\omega)=\sum_{k=0}^{n(\omega)-1}
  \ket{h_k(\omega)}\bra{h_k(\omega)}\in\bddlin^1_+(H),\quad
  \tr\,T(\omega)=1,\quad\rank(T(\omega))=n(\omega),
\]
and all functions of $\omega$ are B\^ochner $\mu$-measurable in
their natural range spaces.
\end{proposition}

\begin{proof}
While the result is a specialization of what we showed for general
PSF measures, we provide another proof, which is considerably
shortened by the use of the well-established Hilbert space
operator theory.

Let us define $\Lambda\in\bddlin^1_+(H)$ by (\ref{eq:Lambda}), and
consider the positive trace class operator valued measure
$F(X):=\Lambda^{1/2}E_0(X)\Lambda^{1/2}$. The total variation of
$F$ is
\[
  \sup\sum_{i=1}^n\|F(X_i)\|_{\bddlin^1(H)}
  =\sup\sum_{i=1}^n\tr\,F(X_i)
  =\tr\,F(\Omega)\leq\tr\,\Lambda\,\cdot\|E_0(\Omega)\|,
\]
where the supremum is over all finite partitions
$\Omega=\bigcup_{i=1}^n X_i$, and we made use of the positivity of
$F(X_i)$ and $\Lambda$, and basic properties of the trace.

Thus $F$ is of bounded total variation. It also has the same
null-sets as $E_0$, which in turn has the same null-sets as the
measure $\mu$ constructed in (\ref{eq:def-mu}). Since
$\bddlin^1(H)$, as a separable (for a reference and an easy direct
proof see e.g. ~\cite[p.\ 794]{KLY}) dual space has the
Radon--Nikod\'ym property (e.g.~\cite{DU}, Theorem~III.3.1), we
may apply the vector-valued Radon--Nikod\'ym theorem to deduce the
existence of an $\bddlin^1(H)$-valued B\^ochner-integrable density
$\omega\mapsto T(\omega)$ such that $F(X)=\int_X
T(\omega)d\mu(\omega)$ for all $X\in\Sigma$. Since $F(X)$ is a
positive operator, also $T(\omega)$ must be for a.e.\
$\omega\in\Omega$; this follows in our separable situation from
the corresponding result for scalar-valued measures and densities,
since the positivity of an operator $A\in\bddlin(H)$ can be tested
in terms of the positivity of $\<\phi|A\phi\>$, where $\phi$ goes
through a countable dense subset of $H$. Moreover,
\[
  \int_X\tr\,T(\omega)\,d\mu(\omega)
  =\tr\,F(X)
  =\sum_{n=1}^{\infty}\frac{\alpha_n}{
    1+\<e_n|E_0(\Omega)e_n\>}\<e_n|E_0(X)e_n\>
  =\mu(X),
\]
and this implies that $\tr\,T(\omega)=1$ for a.e.\
$\omega\in\Omega$.

What remains is the decomposition of $T(\omega)$. For each single
operator, this of course is well known from the Hilbert space
operator theory, but the point is now to have this in a measurable
way. Such a measurable decomposition is proved in \cite{DPZ},
Proposition~1.8, based on a classical theorem on \emph{measurable
selectors} \cite{KRN}. The result in \cite{DPZ} also gives the
additional properties $\|h_{k-1}(\omega)\|\geq\|h_k(\omega)\|>0$
for $1\leq k<n(\omega)$ and $\<h_k(\omega)|h_{\ell}(\omega)\>=0$
for $k\neq\ell$, in addition to those stated in the assertions.
\end{proof}

\begin{remark} Part of the information in the above Proposition
was obtained in a different way in \cite{Beukema}, Proposition 27,
Remark 28, and the proof of Theorem 79.
\end{remark}

\section{Generalized eigenvectors of normal operators}\label{sec:normal-ops}

Here we make a further specialization of the general theory to the
case of a spectral measure associated to a normal operator $T$ in
a Hilbert space. We are going to show that in this situation the
functionals $\ket{d_k(\omega)}$ may be interpreted as generalized
eigenvectors of $T$, in a precise sense to be defined below.
Results of this kind have a long history; instead of attempting a
comprehensive record, we would only like to mention the general
setting for much of the early theory provided by C.~Foia\c{s}
\cite{Foias} and the more recent approach due to S.~J.~L.\ van
Eijndhoven and J.~de Graaf \cite{EG}, from which we borrow an
auxiliary result. The technique based on the use of the measure
$\mu$ and the Radon--Nikod\'ym theorem in our approach bears a
certain resemblance to some ideas already present in \cite{Foias}.
Operator densities in connection with positive operator measures
are also used in \cite{Beukema}, but the generalized eigenvalue
problem is not treated there.

Now, let $H$ be a Hilbert space and $V$ a dense subspace of $H$.
Every $\phi\in H$ determines a continuous antilinear functional
$\psi\mapsto\<\psi|\phi\>$. We denote this functional, and also
its restriction to $V$, by $\ket{\phi}$. Thus $\ket{\phi}$ belongs
to $V^\times$, the linear space of all antilinear 
functionals on $V$, and the mapping $\phi\mapsto\ket{\phi}$ is a
linear injection from $H$ into $V^\times$. We often write simply
$V\subset H\subset V^\times$.

In the following definition we assume that $T:\domain(T)\to H$ is
a densely defined linear map in $H$ and let $T^*:\domain(T^*)\to
H$ be its adjoint. We assume that $V\subset \domain(T^*)$ and
$T^*(V)\subset V$. Let us denote by $\tilde T:V^\times\to
V^\times$ the linear map defined by
\begin{equation}\label{eq:tildeT}
  (\tilde TF)(\psi):=F(T^*\psi),\qquad
  F\in V^\times,\ \psi\in V.
\end{equation}
With $F=\ket{\phi}\in\domain(T)$ and $\psi\in V$, this yields
$
  (\tilde T\ket{\phi})(\psi)=\<T^*\psi|\phi\>=\<\psi|T\phi\>
$
Thus $\tilde T$ may be regarded as an extension of $T$ under the
interpretation $V\subset H\subset V^\times$.

\begin{definition}
If $F\in V^\times\setminus\{0\}$ and $\lambda\in\C$ satisfy
$
  \tilde T F=\lambda F,
$
then $F$ is called a {\em generalized eigenvector} of $T$
belonging to the {\em generalized eigenvalue} $\lambda$ of $T$
(relative to $V$).
\end{definition}

The discussion preceding the definition justifies this
terminology: any eigenvalue of $T$ is a generalized eigenvalue.

Let then $T:\domain(T)\to H$ be a normal operator. According to
the well-known \emph{spectral theorem} (see e.g.\ ~\cite{Rudin}),
associated to $T$ there is a uniquely determined spectral measure
$E_0:\mathcal{B}(\C)\to\bddlin(H)$, supported on the spectrum of
$T$, such that
\begin{equation}\label{eq:spectralThm}
  \<\phi|T\psi\>=\int_{\C}\lambda\<\phi|E_0(d\lambda)\psi\>,
  \qquad\phi\in H,\,\psi\in\domain(T).
\end{equation}
The following result is based on (\ref{eq:spectralThm}) and the
application of our diagonalization results to $E_0$. The existence
of invariant subspaces $V$ for $T$ and $T^*$, as required in the
theorem, is always guaranteed by results at the end of the
section.

\begin{theorem}\label{thm:eigen}
Let $T:\domain(T)\to H$ be a normal operator as in
(\ref{eq:spectralThm}) with $V\subset\domain(T)\cap\domain(T^*)$
and $TV\subset V$, $T^*V\subset V$. Then there exist a finite
positive Borel measure $\mu$ on $\sigma(T)$ having the same
null-sets as $E_0$, and a sequence of $\mu$-measurable mappings
$\lambda\mapsto\ket{d_k(\lambda)}\in V^\times$ such that, for
$\mu$-almost every $\lambda\in\sigma(T)$ and every $k\in\N$,
$\ket{d_k(\lambda)}$ is either zero or a simultaneous generalized
eigenvector of $T$ and $T^*$ relative to $V$ belonging to their
generalized eigenvalues $\lambda$ and $\bar{\lambda}$,
respectively. Moreover,
\begin{equation}\label{eq:phiNpsi}
  \<\phi|T\psi\>=\int_{\sigma(T)}\lambda
  \sum_{k=0}^{n(\lambda)-1}
  \<\phi|d_k(\lambda)\>\<d_k(\lambda)|\psi\>d\mu(\lambda),
  \qquad\phi,\ \psi\in V.
\end{equation}
\end{theorem}

\begin{proof}
From (\ref{eq:spectralThm}) we have 
$
  \<\phi|T\psi\>=\<T^*\phi|\psi\>=\int_{\C}\lambda\<\phi|E_0(d\lambda)\psi\>,$
$\phi,\,\psi\in V.$
On the other hand, we know the existence of $\mu$ and
$\ket{d_k(\lambda)}$ so that
$
  \<\phi|E_0(d\lambda)\psi\>=\sum_{k=0}^{n(\lambda)-1}
  \<\phi|d_k(\lambda)\>\<d_k(\lambda)|\psi\>d\mu(\lambda),$
  $\phi,\,\psi\in V,$
and hence
\begin{equation}\label{eq:phiNpsi}
  \<\phi|T\psi\>=\<T^*\phi|\psi\>=\int_{\C}\lambda
  \sum_{k=0}^{n(\lambda)-1}
  \<\phi|d_k(\lambda)\>\<d_k(\lambda)|\psi\>d\mu(\lambda),
  \qquad\phi,\,\psi\in V.
\end{equation}
Since $E_0(\C)=I$, we also have
\begin{equation}\label{eq:phipsi}
  \<\phi|\psi\>=\int_{\C}\sum_{k=0}^{n(\lambda)-1}
  \<\phi|d_k(\lambda)\>\<d_k(\lambda)|\psi\>d\mu(\lambda),
  \qquad\phi,\,\psi\in V.
\end{equation}
Comparing (\ref{eq:phiNpsi}) and (\ref{eq:phipsi}) with $\phi$
replaced by $T^*\phi$ (still in $V$ by the assumption that
$T^*V\subset V$), and recalling that the functions
$\lambda\mapsto(\<d_k(\lambda)|\psi\>)_{k=0}^{\infty}$, $\psi\in
V$, are dense in $L^2_{n(\cdot)}(\Omega,\mu,\ell^2)$, we deduce
that there must hold, for $\mu$-a.e. $\lambda\in\C$,
\begin{equation}\label{eq:eigenvector}
  \<T^*\phi|d_k(\lambda)\>=\lambda\<\phi|d_k(\lambda)\>,
  \qquad\phi\in V.
\end{equation}
By choosing a null-set $N$ such that (\ref{eq:eigenvector}) holds
in its complement for all $\phi=e_n$, $n\in\N$, it also holds, by
linearity, for all $\phi\in V$. Thus, for $\mu$-a.e.
$\lambda\in\sigma(T)$, we have the generalized eigenvector
equations
$
  \tilde{T}\ket{d_k(\lambda)}=\lambda\ket{d_k(\lambda)}.
$
Making a similar comparison of (\ref{eq:phiNpsi}) and
(\ref{eq:phipsi}) with $\psi$ replaced by $T\psi$ (where the
assumption $TV\subset V$ appears), and proceeding analogously, we
deduce the adjoint equation
$
  \widetilde{T^*}\ket{d_k(\lambda)}=\bar{\lambda}\ket{d_k(\lambda)},
$
which completes the proof.
\end{proof}

\begin{remark}\label{rem:sizeOfEigen}
In the construction of the measure $\mu$ in (\ref{eq:def-mu}), any
positive sequence $(\alpha_n)_{n=0}^{\infty}\in\ell^1$ could be
chosen. For a spectral measure $E_0$, we have
$\<e_n|E_0(\C)e_n\>=1$ for all $n\in\N$, so that the eigenvalues
$\beta_n$ of the operator $\Lambda\in\bddlin^1_+(H)$ in
(\ref{eq:Lambda}) are just $\beta_n=\alpha_n/2$. Since we know
that the action of the $\ket{d_k(\lambda)}$ can be extended
boundedly to $\range(\Lambda^{1/2})$, and since
$(\beta_n)_{n=0}^{\infty}\mapsto(\beta_n^{1/2})_{n=0}^{\infty}$ is
a bijection between the positive cones of $\ell^1$ and $\ell^2$,
we have the following information concerning the size of the
generalized eigenvectors (cf.~(\ref{eq:sizeOfDk})): For any
positive $(\varrho_n)_{n=0}^{\infty}\in\ell^2$, for a.e.\
$\lambda\in\sigma(T)$ (with respect to the spectral measure of the
normal operator $T$), there exist $\ket{d_k(\lambda)}\in
V^{\times}$, $0\leq k<n(\lambda)$, which are simultaneous
generalized eigenvectors of $T$ and $T^*$, and satisfy
$
  \big(\varrho_n\<e_n|d_k(\lambda)\>\big)_{n=0}^{\infty}
  \in\ell^2.
$
\end{remark}

The generalized eigenvector equations can be used, at least in
principle, to solve for the $d_k(\lambda)$, from which one may try
to construct the spectral decomposition of $T$. We illustrate this
in a toy example below, but also point out in the following
section some intrinsic problems related to this approach.

\begin{example}\label{shiftOperator}
Consider, as in Example~\ref{ex:shifts}, a Hilbert space $H$ with
an orthonormal basis $(e_n)_{n\in\Z}$ and the simplest shift
operator $S:e_n\mapsto e_{n-1}$ (hence $S^*:e_n\mapsto e_{n+1}$).
Clearly $S$ is unitary, in particular normal, so that the theory
developed in this section applies to it.
 Writing
$\ket{d}=\sum_{j=-\infty}^{\infty}d_j\ket{e_j}$ for
$\ket{d_k(\lambda)}$, (\ref{eq:eigenvector}) reads for $\phi=e_n$
as $d_{n+1}=\lambda d_n$, giving the unique (up to normalization)
solution $d_j=d_0\lambda^j$, except for the case $\lambda=0$,
where $d_j\equiv 0$ is the only solution. Thus, for every
$\lambda\in\C\setminus\{0\}$, there corresponds a one-dimensional
generalized eigenspace spanned by
$\ket{d(\lambda)}=\sum_{j=-\infty}^{\infty}\lambda^{j}\ket{e_j}$.

The eigenvector equation for the adjoint $S^*$ gives similarly
$d_{n-1}=\bar{\lambda} d_n$, yielding the solution
$d_j=\bar{\lambda}^{-j} d_0$. This can only coincide with the
generalized eigenvector of $S$ if $\bar{\lambda}=\lambda^{-1}$.
Thus the only simultaneous generalized eigenvectors of $S$ and
$S^*$ are, up to normalization,
$\ket{d(e^{it})}=\sum_{j=-\infty}^{\infty}e^{itj}\ket{e_j}$. Hence
the general theory guarantees that, for some finite positive
measure $\mu$ on the unit circle $\T$, the spectral measure of $S$
is given by
$
  \<\phi|E_S(X)\psi\>
  =\int_{X\cap\T}
  \<\phi|d(\lambda)\>\<d(\lambda)|\psi\>d\mu(\lambda)$
  for all
  $\phi,\ \psi\in V.$
Testing the equality $\<\phi|E_S(\T)\psi\>=\<\phi|\psi\>$ with
$\phi=e_m$, $\psi=e_n$, we find that
$
  \int_{\T}\lambda^{m-n}d\mu(\lambda)=\delta_{m,n}$ 
for all  $m,\,n\in\T.$
This shows that the Fourier coefficients of the measure $\mu$
coincide with those of the normalized Haar measure $\nu$ of $\T$,
and hence
$
  \<\phi|E_S(X)\psi\>=\int_{X\cap\T}
  \<\phi|d(\lambda)\>\<d(\lambda)|\psi\>d\nu(\lambda)
  =\sum_{m,n\in\Z}\int_{X\cap\T} \lambda^{m-n}d\nu(\lambda)
  \<\phi|e_m\>\<e_n|\psi\>$ for all
  $\phi,\ \psi\in V,$
gives the spectral measure of $S$.
\end{example}
We finally address the question of invariant subspaces as required
in Theorem~\ref{thm:eigen}. Following \cite{EG}, we define the
\emph{joint $\Cinfty$-domain} of linear operators
$A_i:\domain(A_i)\subset H\to H$, $i=1,\ldots,k$, as
$\Cinfty(A_1,\ldots,A_k):=\big\{\psi\in H|\text{ for all }N\in\N\text{ and }
     \pi\in\{1,\ldots,k\}^N,\
     \psi\in\domain(A_{\pi(1)}A_{\pi(2)}\cdots A_{\pi(N)})\big\}.$
We quote the following result and deduce an immediate consequence:

\begin{theorem}[\cite{EG}]
Suppose that $\Cinfty(A_1,\ldots,A_k)$ is dense in $H$. Then there
exists an orthonormal basis $(e_n)_{n=0}^{\infty}$ such that each
$A_i$, $i=1,\ldots,k$, maps $V=\lin\{e_n|n\in\N\}$ into itself.
\end{theorem}

\begin{corollary}
Let $T:\domain(T)\subset H\to H$ be normal \textup{(}or
bounded\textup{)}. Then there exists an orthonormal basis
$(e_{n})_{n=0}^{\infty}$ such that both $T$ and $T^*$ map $V$ into
itself.
\end{corollary}

\begin{proof}
If $T$ is normal, let $E_0$ be its spectral measure and $D_k$ the
disc $\{\zeta\in\C||\zeta|<k\}$, every $\phi\in H$ satisfies
$E_0(D_k)\phi\in\Cinfty(T,T^*)$ for all $k\in\N$ and
$E_0(D_k)\phi\to\phi$ as $k\to\infty$. Hence $\Cinfty(T,T^*)$ is
dense in $H$, and we can apply the previous theorem.
If $T$ is bounded, then $\Cinfty(T,T^*)=H$, and we derive the same
conclusion.
\end{proof}

\section{A counterexample concerning generalized
eigenvectors}\label{sec:counterexample}

In Theorem~\ref{thm:eigen} we proved that almost every (with
respect to the spectral measure) point $\lambda$ in the spectrum
of a normal operator $T$ is a generalized eigenvalue of $T$, and
moreover Remark~\ref{rem:sizeOfEigen} showed that the associated
generalized eigenvectors are in a sense not very far from being
vectors in the Hilbert space $H$. In this section we show that the
converse statement fails: even if some $\lambda\in\C$ is a
generalized eigenvalue of $T$ with a ``nice'' associated
generalized eigenvector, this $\lambda$ need not be in the Hilbert
space spectrum of $T$, and this can already happen for a bounded
(in fact, Hilbert--Schmidt) self-adjoint operator $T$. The
following technical lemma is the key to the counterexample:

\begin{lemma}
There exists an infinite matrix $(m_{ij})_{i,j=0}^{\infty}$ having
entries in $\left[0,1\right]$ and with the following properties:
\begin{itemize}
  \item $m_{ij}=m_{ji}$ for all $i,j\in\N$,
  \item for all $i\in\N$, the sequence $(m_{ij})_{j=0}^{\infty}$
  has only finitely many non-zero members,
  \item for all $i\in\N$, there holds
  $\sum_{j=0}^{\infty}m_{ij}=1$, and
  \item $\sum_{i,j=0}^{\infty}m_{ij}^2<1$.
\end{itemize}
\end{lemma}

\begin{proof}
We start with a sequence $(a_i)_{i=0}^{\infty}$ of positive
integers, to be chosen later. Denote $\sigma_{-1}:=0$ and
$\sigma_i:=\sum_{k=0}^i a_k$ for $i\in\N$. Let $\phi(j)$, for
$j=1,2,\ldots$, be the unique $i\in\N$ such that
$\sigma_{i-1}<j\leq\sigma_i$. Since $\sigma_{i-1}\geq i$ (as
$a_k\geq 1$), we have $\phi(j)<j$. Then we define a sequence
$(b_j)_{j=0}^{\infty}\in[0,1]^{\N}$ inductively as follows:
$b_0:=1/a_0$, and $b_j:=(1-b_{\phi(j)})/a_j$ for $j=1,2,\ldots$.
Finally, the matrix entries $m_{ij}$ are defined: We set
$m_{ii}:=0$, and for $j>i$ we let $m_{ij}:=b_i$ if
$\sigma_{i-1}<j\leq\sigma_i$ (i.e., $\phi(j)=i$), and zero
otherwise. The entries $m_{ij}$ with $j<i$ are defined so as to
satisfy the symmetry requirement.

That $(m_{ij})_{j=0}^{\infty}$ has only finitely many non-zero
members is clear from this definition. In fact, for $j\geq i$,
there are $\sigma_i-\sigma_{i-1}=a_i$ entries equal to
$b_i=(1-b_{\phi(i)})/a_i$, the others being zero. For $j<i$, an
entry $m_{ij}=m_{ji}$ can only differ from zero if $j=\phi(i)$, in
which case it is equal to $b_j=b_{\phi(i)}$. Thus
$
  \sum_{j=0}^{\infty}m_{ij}
  =b_{\phi(i)}+a_i\times(1-b_{\phi(i)})/a_i=1,
$
as we wanted. To compute the square sum of the entries $m_{ij}$,
observe that there are exactly $2(\sigma_i-\sigma_{i-1})=2a_i$
entries for which we gave the value $b_i$; hence
$
  \sum_{i,j=0}^{\infty}m_{ij}^2
  =\sum_{i=0}^{\infty}2a_i b_i^2
  =\sum_{i=0}^{\infty}2a_i\big(1-b_{\phi(i)}\big)^2/a_i^2
  \leq\sum_{i=0}^{\infty}2/a_i,
$
and this is easily made to satisfy the final requirement by a
suitable choice of $(a_i)_{i=0}^{\infty}$.
\end{proof}

We now provide the announced counterexample. Once again, let $H$
be a Hilbert space with an orthonormal basis
$(e_n)_{n=0}^{\infty}$ and $V:=\lin(e_n)_{n=0}^{\infty}$. We
define an operator $T\in\bddlin(H)$ in terms of the matrix
$(m_{ij})_{i,j=0}^{\infty}$, i.e., we set
\begin{equation}\label{def:TbyMatrix}
  T\Big(\sum_{i=0}^{\infty}a_i e_i\Big)
  :=\sum_{i=0}^{\infty}\Big(\sum_{j=0}^{\infty}m_{ij}a_j\Big)
    e_i
\end{equation}
Since the matrix is real and symmetric, the operator $T$ is
self-adjoint, and hence $\sigma(T)\subset\R$. Moreover, $T$ is a Hilbert--Schmidt operator, and
$
  \|T\|_{\bddlin^2(H)}^2=\sum_{\lambda\in\sigma(T)}\lambda^2
  =\sum_{i,j=0}^{\infty}m_{ij}^2<1.
$
In particular, $T$ has a discrete spectrum, and its eigenvalues
$\lambda$ satisfy
$
  -1<\lambda_{-}\leq\lambda\leq\lambda_{+}<1.
$

However, let us consider the generalized eigenvalue problem for
$T$. The fact that every column and row of the infinite matrix
 $(m_{ij})_{i,j=0}^{\infty}$ has only finitely many non-zero entries is
 equivalent with the properties $TV\subset V$ and $T^* V\subset V$. Thus
 the generalized eigenvector formalism of the previous section is applicable.

Observe that the same extension $\tilde{T}$ on $V^\times$ is
obtained by using the original defining
formula~(\ref{def:TbyMatrix}), which also makes sense for an
arbitrary $\sum_{j=0}^{\infty}a_j\ket{e_j}\in V^\times$. We may
use this remark to compute $\tilde{T}\ket{e}$, where
$\ket{e}:=\sum_{j=0}^{\infty}\ket{e_j}$. From self-adjointness and
the fact that $\sum_{j=0}^{\infty}m_{ij}=1$, it follows at once
that $\tilde{T}\ket{e}=\widetilde{T^*}\ket{e}=\ket{e}$, and hence
$\ket{e}$ is a generalized eigenvector of $T$ and $T^*$
corresponding to the generalized eigenvalue $1$. By the remarks
concerning the spectrum of $T$ above, this generalized eigenvalue
is not in $\sigma(T)$. Summarizing, we have shown the following:

\begin{proposition}
There exists a self-adjoint operator $T\in\bddlin^2(H)$ such that
$TV=T^*V\subset V$, and hence extending to
$\tilde{T}:V^{\times}\to V^{\times}$ by (\ref{eq:tildeT}), such
that
\begin{itemize}
  \item $\|T\|_{\bddlin^2(H)}<1$; in particular,
$\sigma(T)\subset[\lambda_{-},\lambda_{+}]\subset\left]-1,1\right[$,
but
  \item $\tilde{T}\ket{e}=\widetilde{T^*}\ket{e}=\ket{e}$ where
$\ket{e}=\sum_{j=0}^{\infty}\ket{e_j}$.
\end{itemize}
\end{proposition}

\end{document}